\documentclass[12pt]{amsart}

\usepackage{geometry}
\usepackage{amsfonts,amsmath,amssymb,amsthm}
\usepackage{euscript,mathrsfs}
\usepackage{enumerate}
\usepackage[colorlinks,citecolor=red,linkcolor=blue,hyperindex]{hyperref}

\newtheorem{theorem}{Theorem}[section]
\newtheorem{lemma}[theorem]{Lemma}
\newtheorem{proposition}[theorem]{Proposition}
\newtheorem{corollary}[theorem]{Corollary}

\theoremstyle{definition}

\theoremstyle{remark}

\newtheorem*{remark*}{Remark}

\newcommand{\N}{\mathbb{N}}
\newcommand{\Z}{\mathbb{Z}}

\newcommand{\C}{\mathbb{C}}
\newcommand{\K}{\mathcal{K}}

\newcommand{\BH}{\mathcal{B}(H)}
\newcommand{\QH}{\mathcal{Q}(H)}
\newcommand{\Fad}{\mathcal{F}^{\mathrm{ad}}_n(H)}
\newcommand{\Fen}{\mathcal{F}^{\mathrm{en}}_n(H)}

\DeclareMathOperator{\Ran}{Ran}
\DeclareMathOperator{\Ker}{Ker}
\newcommand{\spann}{\mathrm{span}}

\DeclareMathOperator{\ind}{ind}

\numberwithin{equation}{section}

\begin{document}

\title[Almost doubly commuting Fredholm tuples]
{Path Components of Almost Doubly Commuting Fredholm Tuples}

\author{Chen Li}
\address{Chen Li: School of Mathematical Sciences, Fudan University,
Shanghai, China}
\email{22110180024@m.fudan.edu.cn}

\author{Guohao Li}
\address{Guohao Li: School of Mathematical Sciences, Fudan University,
Shanghai, China}
\email{24110180024@m.fudan.edu.cn}

\author{Kai Wang}
\address{Kai Wang: School of Mathematical Sciences, Fudan University,
Shanghai, China}
\email{kwang@fudan.edu.cn}
\date{November 10, 2026}

\begin{abstract}
In this article, we prove that the Fredholm index classifies the
norm-path components of almost doubly commuting Fredholm $n$-tuples
on a separable infinite-dimensional Hilbert space. Using Kasparov's
technical theorem, we join each tuple to an essentially normal one
by a straight-line Fredholm path and apply the Douglas--Curto theorem.
\end{abstract}


\maketitle


\section{Introduction}\label{sec-introduction}

Let $H$ be a separable infinite-dimensional complex
Hilbert space, and write $\BH$ and $\K(H)$ for the bounded and
compact operators on $H$, respectively. An operator $T\in\BH$ is
Fredholm if its range is closed and both $\Ker T$ and $\Ker T^*$
are finite dimensional. Its index is
\[
 \ind T=\dim\Ker T-\dim\Ker T^*.
\]

Compact perturbations preserve both Fredholmness and the index.
Atkinson's theorem \cite{Atkinson} gives an algebraic characterization
of Fredholmness: $T$ is Fredholm if and only if its image is invertible
in the Calkin algebra
\[
 \QH=\BH/\K(H),\qquad q:\BH\longrightarrow\QH.
\]

Applied to $T-\lambda I$, this criterion defines the index function
$\lambda\mapsto\ind(T-\lambda I)$ outside the essential spectrum
$\sigma_{\mathrm e}(T)=\sigma(q(T))$.
The Brown--Douglas--Fillmore theorem \cite{BDF1973} states that
two essentially normal operators $S,T\in\BH$ are essentially
unitarily equivalent if and only if they have the same essential
spectrum and the same index function on its complement.

The classification by norm-continuous paths uses the Fredholm
index alone. The Fredholm operators form an open subset of $\BH$,
on which the index is locally constant. Conversely, two Fredholm
operators on $H$ with the same index can be joined by a
norm-continuous Fredholm path \cite{Atiyah}.
Since every integer occurs as an index, their path components
are indexed by $\Z$.

For tuples, the Fredholm deformation problem requires a notion of
joint invertibility and paths that preserve the commutation
relations between the coordinates.
Using Taylor's joint spectrum \cite{Taylor},
Curto \cite{Curto} developed the corresponding Fredholm theory
for operator tuples.
Fix $n\geq1$. A tuple $A=(A_1,\ldots,A_n)$ is
\emph{almost commuting} if $[A_i,A_j]\in\K(H)$ for all $i,j$.
Its image $a=(a_1,\ldots,a_n)$, where $a_i=q(A_i)$, is then a
commuting tuple in $\QH$. Such a tuple $A$ is \emph{Fredholm}
if $0\notin\sigma_T(a)$, where $\sigma_T$ denotes the Taylor spectrum.

The tuple $A$ is \emph{almost doubly commuting} if
\begin{align}\label{def-almost-double}
 [A_i,A_j],\ [A_i,A_j^*]\in\K(H)\qquad(i\ne j),
\end{align}
and it is \emph{essentially normal} if, in addition,
\begin{align}\label{def-essential-normal}
 [A_i,A_i^*]\in\K(H)\qquad(1\leq i\leq n).
\end{align}
We write $\Fad$ for the almost doubly commuting Fredholm
$n$-tuples and $\Fen$ for the essentially normal Fredholm
$n$-tuples. 

Curto's study of connectedness for Taylor-invertible tuples
\cite{CurtoConnected} was followed by his Fredholm theory for
almost commuting tuples \cite{Curto}. He obtained a multivariable
Atkinson theorem and defined an integer-valued index that is
invariant under compact perturbations and norm-continuous Fredholm
paths, and additive under orthogonal direct sums. He also computed
the indices of the coordinate tuples on the Hardy spaces of the
ball and the polydisc. We follow his index convention, recalled
in Section~\ref{sec-koszul}.

Curto's related spectral work includes spectral inclusion for
doubly commuting subnormal tuples \cite{CurtoSubnormal}, spectral
permanence in $C^*$-algebras \cite{CurtoPermanence}, and, with
Fialkow, spectral and index formulas for left and right multiplication
tuples \cite{CurtoFialkow}. His survey \cite{CurtoSurvey} gives an
account of Taylor and Harte spectra. His compact-perturbation
question for commuting Fredholm pairs of index zero
\cite{CurtoProblems} received a negative answer from Gelca
\cite{Gelca}. Later work includes a Weyl theorem for commuting
hyponormal pairs under the quasitriangular property
\cite{ChavanCurto} and preservation of Taylor spectra, Taylor
essential spectra and Fredholm indices under the generalized
spherical Aluthge transform of commuting tuples \cite{BCLY}.

The deformation problem was posed by Douglas in
\emph{The relation of Ext to $K$-theory} \cite{Douglas}.
It asks whether two almost commuting Fredholm $n$-tuples with the
same index can be joined by a norm-continuous path of such tuples.
For essentially normal Fredholm tuples, Douglas \cite{Douglas}
proved that two tuples with the same index can be joined by a
norm-continuous path of essentially normal Fredholm tuples.
Curto gave a detailed proof in \cite[Theorem 7]{Curto}.
The proof uses Brown--Douglas--Fillmore extension theory
\cite{BDF1977} applied to the commutative $C^*$-algebra generated
by the Calkin coordinates.

Almost double commutation retains the relations between distinct
coordinates but places no normality condition on an individual
Calkin coordinate. Thus the algebra generated by the Calkin image
need not be commutative, and the essentially normal result does
not directly give the classification in $\Fad$.
For almost doubly commuting Fredholm pairs, Curto
\cite[Theorem 9]{Curto} proved that equal indices suffice
when each endpoint has a semi-Fredholm coordinate. His proof uses
partial isometries and shift models, together with a path between
the coordinate tuples on the Hardy spaces of the ball and the
polydisc \cite{Curto}.

Curto \cite{Curto} noted that the additional semi-Fredholm condition
might be associated with a second invariant. In this article we
show that the index alone suffices in $\Fad$: the coordinate
assumption can be removed, and the classification holds for every
finite $n$.

\begin{theorem}\label{thm-components}
For every $A,B\in\Fad$, the following conditions are equivalent:
\begin{enumerate}[(1)]
\item $\ind A=\ind B$;
\item $A$ and $B$ are joined by a norm-continuous path in $\Fad$.
\end{enumerate}
Every integer occurs as an index. Hence the index induces a bijection
\[
 \pi_0(\Fad)\longrightarrow\Z,
\]
where $\pi_0$ denotes the set of path components.
\end{theorem}

To obtain an essentially normal representative in each component,
we construct a line segment
from any given tuple to $\Fen$. The construction also gives a
uniform lower bound for the Calkin Laplacians.

\begin{theorem}\label{thm-normalization}
For every $A\in\Fad$, there is $C\in\Fen$ such that
\[
 A(t)=(1-t)A+tC\in\Fad\qquad(0\leq t\leq1)
\]
and
\[
 \sum_{i=1}^n|q(C_i)|=1.
\]
Moreover, there is $\eta>0$ such that for every $t\in[0,1]$ and every
choice
\[
 X_i(t)\in\{q(A_i(t))^*q(A_i(t)),\
                  q(A_i(t))q(A_i(t))^*\},
\]
one has $\sum_iX_i(t)\geq\eta1$.
\end{theorem}

Together with the Douglas--Curto theorem, this normalization shows
that the inclusion $\Fen\subseteq\Fad$ induces a bijection of
path components.

Theorem~\ref{thm-normalization} is proved in the Calkin algebra
using a finite separation lemma derived from Kasparov's technical
theorem \cite{Kasparov,Higson}. The lemma supplies commuting
positive weights for a simultaneous normalization of the
coordinates. 

The paper is arranged as follows. Firstly, Section~\ref{sec-koszul} recalls the Koszul complex and Fredholm index theory.
The finite separation lemma will be proved in
Section~\ref{sec-separation}, and the normalization is carried out
in Section~\ref{sec-normalization}. In the end, section~\ref{sec-classification}
completes the proofs of the main theorems.


\section{Koszul operators and the Fredholm index}\label{sec-koszul}
We recall Curto's Fredholm theory \cite{Curto} in the notation
used below.

\subsection{Koszul operators}
Let $e_1,\ldots,e_n$ be the standard orthonormal basis of $\C^n$.
For $I=\{i_1<\cdots<i_k\}$, put
$e_I=e_{i_1}\wedge\cdots\wedge e_{i_k}$ and $e_\varnothing=1$.
Give $\Lambda^\bullet\C^n$ the inner product for which the $e_I$
form an orthonormal basis. Write $\varepsilon_i$ for left exterior
multiplication by $e_i$ and $\iota_i=\varepsilon_i^*$ for contraction:
\[
 \iota_i e_I=
 \begin{cases}
  (-1)^{r-1}e_{I\setminus\{i\}},&i=i_r\in I,\\
  0,&i\notin I.
 \end{cases}
\]
These operators satisfy
\[
 \iota_i\iota_j+\iota_j\iota_i=0,\qquad
 \varepsilon_i\varepsilon_j+\varepsilon_j\varepsilon_i=0,\qquad
 \iota_i\varepsilon_j+\varepsilon_j\iota_i=\delta_{ij}1.
\]
The operators $\varepsilon_i\iota_i$ and $\iota_i\varepsilon_i$
are the orthogonal projections onto the spans of the basis vectors
containing and not containing $i$, respectively.

Set $\mathcal E_k=H\otimes\Lambda^k\C^n$ and
$\mathcal E=\bigoplus_{k=0}^n\mathcal E_k$. For a bounded tuple $A$, put
\[
 d_A=\sum_{i=1}^n A_i\otimes\iota_i,\qquad
 d_A^*=\sum_{i=1}^n A_i^*\otimes\varepsilon_i,\qquad
 D_A=d_A+d_A^*.
\]
The anticommutation relations give
\[
 d_A^2=\sum_{i<j}[A_i,A_j]\otimes\iota_i\iota_j.
\]
For a commuting tuple, $d_A^2=0$, and $d_A$ is the boundary
operator of the Koszul complex
\[
 K(A;H):\quad
 0\longrightarrow\mathcal E_n\xrightarrow{d_{A,n}}\mathcal E_{n-1}
 \longrightarrow\cdots\longrightarrow\mathcal E_1
 \xrightarrow{d_{A,1}}\mathcal E_0\longrightarrow0,
\]
where $d_{A,k}=d_A|_{\mathcal E_k}$. The Taylor spectrum is defined as
\[
 \sigma_T(A)=\{\lambda\in\C^n:K(A-\lambda;H)\text{ is not exact}\},
\]
where $A-\lambda=(A_1-\lambda_11,\ldots,A_n-\lambda_n1)$.
Exactness implies that $\Ran d_{A,k+1}=\Ker d_{A,k}$ in every degree 
and the map is zero at the endpoints. However, if $A$ commutes only modulo
compact operators, $d_A^2$ is compact and the lifted sequence need not be a complex. Thus, we may deal with the problem in the Calkin algebra, which will be discussed in the following section. 

Write $\mathcal E_{\mathrm{even}}$ and $\mathcal E_{\mathrm{odd}}$
for the sums of the even and odd degrees. Since $D_A$ is self-adjoint
and interchanges these subspaces, the operator
\[
 B_A=D_A^{\mathrm{odd}}:
 \mathcal E_{\mathrm{odd}}\longrightarrow\mathcal E_{\mathrm{even}}
\]
gives the block decompositions on
$\mathcal E_{\mathrm{even}}\oplus\mathcal E_{\mathrm{odd}}$:
\[
 D_A=\begin{pmatrix}0&B_A\\B_A^*&0\end{pmatrix},\qquad
 D_A^2=\begin{pmatrix}B_AB_A^*&0\\0&B_A^*B_A\end{pmatrix}.
\]

By spectral permanence \cite{CurtoPermanence}, the Taylor spectrum
of a commuting tuple $b$ in a unital $C^*$-algebra $\mathcal A$ can
be computed in any faithful unital representation. For the
represented tuple, \cite[Proposition 3.4 and Theorem 1]{Curto}
identifies exactness of the Koszul complex with invertibility of
its Koszul operator. Faithfulness therefore gives
\[
 \lambda\notin\sigma_T(b)
 \quad\Longleftrightarrow\quad
 D_{b-\lambda}\text{ is invertible in }M_{2^n}(\mathcal A).
\]

Put $a_i=q(A_i)$. Under the identification
$\mathcal Q(\mathcal E)=M_{2^n}(\QH)$, the images of $d_A$ and $D_A$ are
\[
 d_a=\sum_i a_i\otimes\iota_i,\qquad D_a=d_a+d_a^*.
\]
If $A$ is almost commuting, then $d_a^2=0$. The standard exterior
bases identify each parity subspace with $H^{2^{n-1}}$. We also
write $q$ for the resulting entrywise quotient map on the parity
blocks.
For $\theta\in\{0,1\}^n$, set
\[
 X_i^0=a_i^*a_i,\qquad X_i^1=a_ia_i^*,\qquad
 L_\theta(a)=\sum_{i=1}^nX_i^{\theta_i}.
\]

The following is Curto's criterion for almost doubly commuting
tuples, written in the Calkin algebra \cite[Corollary 3.7]{Curto}.
We include the proof to record the lower bound used below.

\begin{samepage}
\begin{proposition}\label{prop-laplacian}
Suppose that $A$ is almost doubly commuting. The following are equivalent:
\begin{enumerate}
 \item $0\notin\sigma_T(a)$;
 \item $D_A^{\mathrm{odd}}$ is Fredholm;
 \item every $L_\theta(a)$ is invertible in $\QH$;
 \item there is $\gamma>0$ such that
 \begin{equation}\label{eq-koszul-lower-bound}
  L_\theta(a)\geq\gamma1\qquad(\theta\in\{0,1\}^n).
 \end{equation}
\end{enumerate}
\end{proposition}
\end{samepage}

\begin{proof}
By the preceding criterion, condition~(1) is equivalent to
invertibility of $D_a$. Firstly, we give the equivalence of condition ~(1) and ~(3).
Since $d_a^2=(d_a^*)^2=0$, we have
$D_a^2=d_ad_a^*+d_a^*d_a$. For $i\ne j$, the corresponding cross
terms cancel:
\[
 a_ia_j^*\otimes\iota_i\varepsilon_j
 +a_j^*a_i\otimes\varepsilon_j\iota_i
 =[a_i,a_j^*]\otimes\iota_i\varepsilon_j=0.
\]
It follows that
\[
 D_a^2=\sum_{i=1}^n
 \left(a_i^*a_i\otimes\varepsilon_i\iota_i
       +a_ia_i^*\otimes\iota_i\varepsilon_i\right).
\]
On the summand indexed by $I\subseteq\{1,\ldots,n\}$, this is
\[
 \sum_{i\in I}a_i^*a_i+\sum_{i\notin I}a_ia_i^*.
\]
The diagonal blocks of $D_a^2$ are therefore the $2^n$ sums
$L_\theta(a)$, indexed by $\theta\in\{0,1\}^n$. For self-adjoint operator $D_a$, it
is invertible if and only if its square is invertible. Thus, condition ~(1) and ~(3) are equivalent.

Let $b=q(B_A)$. In the parity decomposition, $D_a^2$ has diagonal
blocks $bb^*$ and $b^*b$. If $b$ is invertible, both blocks are
invertible. Conversely, their inverses give a right and a left
inverse for $b$:
\[
 b\bigl[b^*(bb^*)^{-1}\bigr]=1,\qquad
 \bigl[(b^*b)^{-1}b^*\bigr]b=1.
\]
These inverses coincide, so $b$ is invertible. Atkinson's theorem
now proves the equivalence of~(2) and~(3).

Finally, note that each $L_\theta(a)$ is positive. Suppose that $L_\theta(a)$ is invertible for all $\theta\in\{0,1\}^n$ ,
continuous functional calculus gives
\[
 L_\theta(a)\geq\|L_\theta(a)^{-1}\|^{-1}1.
\]
Taking the minimum over $\theta$ gives~(4). Conversely,~(4) implies
$\sigma(L_\theta(a))\subseteq[\gamma,\infty)$ with $\gamma>0$, so every
$L_\theta(a)$ is invertible.
\end{proof}

\subsection{The index and Koszul homology}
For $A\in\Fad$, we use Curto's index \cite[Definition 7.1]{Curto}:
\[
 \ind A=\ind B_A=\dim\Ker B_A-\dim\Ker B_A^*.
\]
It is locally constant on $\Fad$ and invariant under compact
perturbations \cite[Theorem 3]{Curto}. In particular, it is constant
along norm-continuous paths in $\Fad$. It is also additive under
finite orthogonal direct sums and invariant under unitary
identifications.

If $A\in\Fad$ is commuting, its Koszul boundary maps have closed
range and its homology spaces
\[
 H_k(A)=\Ker d_{A,k}/\Ran d_{A,k+1}
\]
are finite dimensional \cite[Theorem 2]{Curto}. In this case,
\cite[Corollary 7.2]{Curto} gives
\begin{equation}\label{eq-koszul-euler}
 \ind A=\sum_{k=0}^n(-1)^{k+1}\dim H_k(A).
\end{equation}
\section{A finite separation lemma in the Calkin algebra}\label{sec-separation}

In this section, all ideals are norm-closed and two-sided. For ideals $I_1,\ldots,I_m$
of a $C^*$-algebra $D$, write
\begin{align}\label{def-ideal-product}
 I_1\cdots I_m
 =\overline{\spann}\{x_1\cdots x_m:x_i\in I_i\}.
\end{align}
For separation in corona algebras and its applications to lifting,
see \cite{Pedersen,OlsenPedersen}. Here the separating element must
also commute with a prescribed separable algebra. We use the
following form of Kasparov's technical theorem
\cite{Kasparov,Higson,Blackadar}.
A later proof is given in
\cite[Theorem 1.4 and Proposition 2.8]{FarahHart}.

\begin{theorem}\label{thm-two-ideals}
Let $D\subseteq\QH$ be a separable unital $C^*$-subalgebra, and let
$I,J$ be ideals of $D$ with $IJ=0$. There exists a positive contraction
$e\in D'\cap\QH$ such that
\[
 ex=x\quad(x\in I),\qquad ey=0\quad(y\in J).
\]
\end{theorem}

\begin{lemma}\label{lem-partition}
Let $D\subseteq\QH$ be a separable unital $C^*$-subalgebra, and suppose
that $I_1,\ldots,I_m$ are ideals of $D$ with $I_1\cdots I_m=0$.
There are mutually commuting positive contractions
$\alpha_1,\ldots,\alpha_m\in D'\cap\QH$ such that
\begin{align}\label{eq-ideal-partition}
 \sum_{i=1}^m\alpha_i=1,\qquad \alpha_iI_i=0\quad(1\leq i\leq m).
\end{align}
\end{lemma}

\begin{proof}
We use induction on $m$. Suppose $m\geq2$ and put
$J=I_2\cdots I_m$, a closed ideal of $D$.
For $x\in I_1$, the hypothesis gives $xx_2\cdots x_m=0$ whenever
$x_i\in I_i$. By linearity and norm continuity, $xJ=0$.
Thus $I_1J=0$, and Theorem~\ref{thm-two-ideals} gives a positive
contraction $e\in D'\cap\QH$ such that
\begin{align}\label{eq-induction-separator}
 ex=x\quad(x\in I_1),\qquad eJ=0.
\end{align}

Set $D_1=C^*(D,e)$. It is separable and unital, and $e$ is central
in $D_1$. For $i\geq2$, let $J_i$ be the ideal of $D_1$ generated
by $eI_i$. We claim that
\begin{align}\label{eq-generated-ideal-induction}
 J_i=\overline{\spann}\{e^kx:k\geq1,\ x\in I_i\}.
\end{align}
To see this, denote the right-hand side by $M_i$.
Since $e$ commutes with $D$, one has $(e^kx)^*=e^kx^*$ and,
for $a\in D$,
\[
 a(e^kx)=e^k(ax),\qquad (e^kx)a=e^k(xa).
\]
Here $x^*,ax,xa\in I_i$. Multiplication by $e$ also preserves
$M_i$, since it increases the exponent by one. Thus $M_i$ is a
closed ideal of $D_1$ containing $eI_i$, which derives $J_i\subseteq M_i$.
Conversely, $e^kx=e^{k-1}(ex)\in J_i$ and gives the reverse
inclusion.

For $x_i\in I_i$ and $k_i\geq1$, commutation with $e$ shows that
\[
 (e^{k_2}x_2)\cdots(e^{k_m}x_m)
 =e^{k_2+\cdots+k_m}x_2\cdots x_m=0,
\]
since $x_2\cdots x_m\in J$ and $eJ=0$.
Expanding finite sums and then passing to norm limits in each
factor, we have $J_2\cdots J_m=0$. The induction hypothesis in
$D_1$ therefore gives mutually commuting positive contractions
\[
 \beta_2,\ldots,\beta_m\in D_1'\cap\QH,\qquad
 \sum_{i=2}^m\beta_i=1,\qquad\beta_iJ_i=0.
\]
Define
\[
 \alpha_1=1-e,\qquad \alpha_i=e\beta_i\quad(i\geq2).
\]
Each $\beta_i$ commutes with $e$, and hence
\[
 0\leq\alpha_i=e^{1/2}\beta_i e^{1/2}\leq e\leq1
 \qquad(i\geq2).
\]
By definition, each $\alpha_i$ belong to $D'\cap\QH$.
They commute with $\alpha_1$, and for $i,j\geq2$,
$[\alpha_i,\alpha_j]=e^2[\beta_i,\beta_j]=0$.
Their sum is
\[
 \sum_{i=1}^m\alpha_i=1-e+e\sum_{i=2}^m\beta_i=1.
\]
Finally,~\eqref{eq-induction-separator} gives $\alpha_1I_1=0$.
If $i\geq2$ and $x\in I_i$, we have $ex\in J_i$, thus
\[
 \alpha_ix=\beta_i(ex)=0.
 \qedhere
\]
\end{proof}


\section{Normalization in the Calkin algebra}\label{sec-normalization}

Let $a=(a_1,\ldots,a_n)$ be a doubly commuting Taylor-invertible tuple
in $\QH$, and put
\[
D_i=C^*(1,a_i),\qquad D=C^*(1,a_1,\ldots,a_n).
\]
The algebra $D$ is separable. Since $a_i$ and $a_i^*$ commute with
$a_j$ and $a_j^*$ for $i\ne j$, the algebras $D_i$ and $D_j$ commute
elementwise. Thus
\begin{align}\label{def-coordinate-algebras}
D=\overline{\spann}\{d_1\cdots d_n:d_i\in D_i\},
\end{align}
since the linear span on the right is a unital $*$-algebra containing
all the generators of $D$.

\subsection{Spectral cutoffs and ideals}
The lower bound in Proposition~\ref{prop-laplacian} is a joint
condition; the individual coordinates need not be invertible.
We therefore use spectral cutoffs for both $a_i^*a_i$ and $a_ia_i^*$
before applying the separation lemma.

Choose $\gamma>0$ as in Proposition~\ref{prop-laplacian}, and fix
\begin{align}\label{def-cutoff}
0<\delta<\frac{\gamma}{n},\qquad
f(s)=(\delta-s)_+,\qquad
g(s)=s+f(s)=\max\{s,\delta\}\quad(s\geq0),
\end{align}
where $r_+=\max\{r,0\}$. Define
\begin{align}\label{def-defect-ideals}
I_i=\operatorname{Ideal}_D\{f(a_i^*a_i),f(a_ia_i^*)\}.
\end{align}
\begin{lemma}\label{lem-low-spectrum}
The ideals in \eqref{def-defect-ideals} satisfy $I_1\cdots I_n=0$.
\end{lemma}

\begin{proof}
Fix $X_i\in\{a_i^*a_i,a_ia_i^*\}$ for each $i$.
These are commuting positive elements. Let $K$ be their joint spectrum
in the commutative algebra $\mathcal A=C^*(1,X_1,\ldots,X_n)$.
Each point of $K$ has the form
$(\chi(X_1),\ldots,\chi(X_n))$ for a character $\chi$ of $\mathcal A$.
Since positivity is unchanged on passing to the unital $C^*$-subalgebra
$\mathcal A$, $X_i\geq0$ and $\sum_iX_i\geq\gamma1$ which implies
$\chi(X_i)\geq0$ and $\sum_i\chi(X_i)\geq\gamma$. Hence
\begin{align}\label{eq-joint-spectrum-region}
K\subseteq
\left\{(s_1,\ldots,s_n)\in[0,\infty)^n:
                       \sum_{i=1}^n s_i\geq\gamma\right\}.
\end{align}
If $\prod_i f(s_i)\ne0$, then $s_i<\delta$ for all $i$, so
$\sum_i s_i<n\delta<\gamma$. Thus $\prod_i f(s_i)$ vanishes on $K$,
and joint continuous functional calculus gives
\begin{align}\label{eq-cutoff-product-zero}
f(X_1)\cdots f(X_n)=0
\end{align}
for all $2^n$ choices of $X_1,\ldots,X_n$.

To pass from these generators to the ideals in
\eqref{def-defect-ideals}, first set
\[
L_i=\operatorname{Ideal}_{D_i}\{f(a_i^*a_i),f(a_ia_i^*)\}.
\]
Since $D_i$ is unital and the generators are self-adjoint, $L_i$
is the norm-closed linear span of the elements
\[
u_i f(X_i)v_i,\qquad
u_i,v_i\in D_i,\quad X_i\in\{a_i^*a_i,a_ia_i^*\}.
\]
Since different $D_i$ commute, \eqref{eq-cutoff-product-zero} gives
\begin{align*}
\prod_{i=1}^n u_i f(X_i)v_i
&=(u_1\cdots u_n)
       \bigl(f(X_1)\cdots f(X_n)\bigr)(v_1\cdots v_n)\\
&=0.
\end{align*}
Only factors from different $D_i$ have been interchanged. By
multilinearity the product is zero when each factor is a finite sum
of these generators. Passing to norm limits gives
\begin{align}\label{eq-coordinate-ideal-product}
\ell_1\cdots\ell_n=0\qquad(\ell_i\in L_i).
\end{align}

We claim that
\begin{align}\label{eq-extended-coordinate-ideal}
I_i=\overline{\spann}
\{d_1\cdots d_{i-1}\ell_i d_{i+1}\cdots d_n:
                    \ell_i\in L_i,\ d_j\in D_j\}.
\end{align}
Let $M_i$ denote the right-hand side. Since different $D_j$ commute,
$M_i$ is self-adjoint. Multiplying a spanning product on either side
by $b\in D_j$ changes only its $j$th factor; when $j=i$, that factor
stays in $L_i$ because $L_i$ is an ideal of $D_i$. Hence $M_i$ is
invariant under multiplication by each $D_j$, and therefore by $D$.
Thus $M_i$ is a closed ideal of $D$.
It contains both generators of $I_i$, so $I_i\subseteq M_i$.
Conversely, $I_i$ contains $L_i$ and is an ideal of $D$; it therefore
contains every product spanning $M_i$. This proves
\eqref{eq-extended-coordinate-ideal}.

For each $i$, take a spanning product
$z_i=z_{i1}\cdots z_{in}$ from \eqref{eq-extended-coordinate-ideal},
with $z_{ij}\in D_j$ and $z_{ii}\in L_i$. Regrouping by coordinate
gives
\[
z_1\cdots z_n=\prod_{j=1}^n(z_{1j}z_{2j}\cdots z_{nj}).
\]
The $j$th parenthesis belongs to $L_j$, since its factors belong to
$D_j$ and $z_{jj}\in L_j$. The product is zero by
\eqref{eq-coordinate-ideal-product}. For arbitrary $x_i\in I_i$,
choose finite linear combinations $x_i^{(r)}$ of these spanning
products such that $x_i^{(r)}\to x_i$ in norm. Each
$x_1^{(r)}\cdots x_n^{(r)}$ is zero, and continuity of multiplication
gives $x_1\cdots x_n=0$. Taking the closed linear span of these
products proves $I_1\cdots I_n=0$.
\end{proof}

From Lemma~\ref{lem-partition}, there exist mutually commuting positive
contractions $\alpha_1,\ldots,\alpha_n\in D'\cap\QH$ satisfying
\begin{align}\label{eq-weight-annihilation}
\sum_{i=1}^n\alpha_i=1,\qquad \alpha_i I_i=0.
\end{align}
In particular,
\begin{align}\label{eq-cutoff-annihilation}
\alpha_i f(a_i^*a_i)=\alpha_i f(a_ia_i^*)=0.
\end{align}
Define
\begin{align}\label{def-normalized-coordinates}
c_i=\alpha_i a_i g(a_i^*a_i)^{-1/2}.
\end{align}

\begin{proposition}\label{prop-endpoint}
The tuple $c=(c_1,\ldots,c_n)$ consists of commuting normal elements
and satisfies
\begin{align}\label{eq-normalized-endpoint}
c_i^*c_i=c_ic_i^*=\alpha_i^2,\qquad
\sum_{i=1}^n|c_i|=1,\qquad
\sum_{i=1}^n c_i^*c_i\geq\frac1n1.
\end{align}
In particular, $c$ is Taylor invertible.
\end{proposition}

\begin{proof}
Firstly, it is known that for any element $a$ of a unital $C^*$-algebra and any continuous
function $h$ on a compact interval containing the spectra of $a^*a$
and $aa^*$, one has
\begin{align}\label{eq-intertwining-functional-calculus}
a h(a^*a)=h(aa^*)a.
\end{align}
Indeed, associativity gives $a(a^*a)^k=(aa^*)^ka$ for every
integer $k\geq0$. The identity follows for polynomials and then for
continuous functions by uniform approximation.

Fix $i$ and write $X=a_i^*a_i$, $Y=a_ia_i^*$ and
$\alpha=\alpha_i$. Since $\alpha$ commutes with $a_i$ and with
the continuous functions of $X$ and $Y$, we obtain
\begin{align*}
c_i^*c_i
&=\alpha^2 g(X)^{-1/2}Xg(X)^{-1/2}\\
&=\alpha^2Xg(X)^{-1}
 =\alpha^2\bigl(1-f(X)g(X)^{-1}\bigr)
 =\alpha^2.
\end{align*}
Here we used $X=g(X)-f(X)$ and
\eqref{eq-cutoff-annihilation}. 
 Choose $h=g^{-1}$ in \eqref{eq-intertwining-functional-calculus}.
 Thus,
\begin{align*}
c_ic_i^*
&=\alpha^2a_i g(X)^{-1}a_i^*
 =\alpha^2g(Y)^{-1}Y\\
&=\alpha^2\bigl(1-f(Y)g(Y)^{-1}\bigr)
 =\alpha^2.
\end{align*}
We get that $c_i$ is normal and $|c_i|=\alpha_i$.

Put $\widetilde D_i=C^*(D_i,\alpha_i)$. 
Hence the algebras $\widetilde D_i$ commute elementwise. For each $a_i$
and $c_i$ belong to $\widetilde D_i$, we have
$[c_i,c_j]=[c_i,c_j^*]=0$ for $i\ne j$.
Consequently, $\sum_i\alpha_i=1$ implies
\begin{align}\label{eq-weight-squares}
0\leq\sum_{i=1}^n\left(\alpha_i-\frac1n1\right)^2
=\sum_{i=1}^n\alpha_i^2-\frac1n1.
\end{align}
Since $c_i^*c_i=c_ic_i^*$, every $L_\theta(c)$ equals
$\sum_i\alpha_i^2$. Therefore, Proposition~\ref{prop-laplacian} gives
Taylor invertibility of the tuple $c$.
\end{proof}

\subsection{The straight-line path}
\begin{proposition}\label{prop-path}
For $0\leq t\leq1$, set $a_i(t)=(1-t)a_i+tc_i$. Then $a(t)$ is
doubly commuting and Taylor invertible. More precisely, put
\begin{align}\label{def-path-constant}
m=\left(\max\{1,\delta,\|a_1\|^2,\ldots,\|a_n\|^2\}\right)^{-1/2}>0.
\end{align}
For every $t\in[0,1]$ and every choice of
$X_i(t)\in\{a_i(t)^*a_i(t),a_i(t)a_i(t)^*\}$, one has
\begin{align}\label{eq-uniform-path-bound}
\sum_{i=1}^n X_i(t)\geq\frac{m^2\delta}{n}1.
\end{align}
\end{proposition}

\begin{proof}
Since $a_i(t)\in\widetilde D_i$, the preceding proof gives
$[a_i(t),a_j(t)]=[a_i(t),a_j(t)^*]=0$ for $i\ne j$.
For either $X_i=a_i^*a_i$ or $X_i=a_ia_i^*$, continuous functional
calculus gives
\[
X_i-\delta1=(X_i-\delta1)_+-f(X_i).
\]
Multiplication by $\alpha_i^2$, followed by
\eqref{eq-cutoff-annihilation}, yields
\begin{align}\label{eq-weighted-spectral-gap}
\alpha_i^2X_i-\delta\alpha_i^2
=\alpha_i^2(X_i-\delta1)_+\geq0,
\end{align}
where positivity follows from the commutativity of the two factors.

To estimate $a_i(t)^*a_i(t)$, write $a_i(t)=a_i k_i(t)$ with
\begin{align}\label{eq-right-factorization}
k_i(t)=(1-t)1+t\alpha_i g(a_i^*a_i)^{-1/2}.
\end{align}
Continuous functional calculus gives that
\[
\|g(a_i^*a_i)\|=\|g(a_ia_i^*)\|
=\max\{\|a_i\|^2,\delta\}\leq m^{-2}.
\]
Thus both $g(a_i^*a_i)^{-1/2}$ and $g(a_ia_i^*)^{-1/2}$ are at
least $m1$. Since $0\leq\alpha_i\leq1$ and $m\leq1$, we have
\[
k_i(t)\geq(1-t)\alpha_i+tm\alpha_i
=\bigl(1-t+tm\bigr)\alpha_i\geq m\alpha_i.
\]
Since $k_i(t)$ and $\alpha_i$ commutes,
\[
k_i(t)^2-m^2\alpha_i^2
=(k_i(t)-m\alpha_i)(k_i(t)+m\alpha_i)\geq0.
\]
Multiplying by the commuting positive element $a_i^*a_i$ and using
\eqref{eq-weighted-spectral-gap}, we obtain
\begin{align}\label{eq-first-orientation-bound}
a_i(t)^*a_i(t)
=a_i^*a_i k_i(t)^2
\geq m^2\alpha_i^2a_i^*a_i
\geq m^2\delta\alpha_i^2.
\end{align}

For $a_i(t)a_i(t)^*$, \eqref{eq-intertwining-functional-calculus} gives
\[
a_i(t)=k'_i(t)a_i,\qquad
k'_i(t)=(1-t)1+t\alpha_i g(a_ia_i^*)^{-1/2}.
\]
 Through the
same argument yields
\begin{align}\label{eq-second-orientation-bound}
a_i(t)a_i(t)^*
=a_ia_i^* k'_i(t)^2
\geq m^2\alpha_i^2a_ia_i^*
\geq m^2\delta\alpha_i^2.
\end{align}
Both choices of $X_i(t)$ therefore satisfy
$X_i(t)\geq m^2\delta\alpha_i^2$. By \eqref{eq-weight-squares},
\[
\sum_{i=1}^nX_i(t)\geq m^2\delta\sum_{i=1}^n\alpha_i^2
\geq\frac{m^2\delta}{n}1.
\]
This proves \eqref{eq-uniform-path-bound}, and Taylor invertibility
follows from Proposition~\ref{prop-laplacian}.
\end{proof}


\section{Proofs of the main theorems}\label{sec-classification}

\begin{proof}[Proof of Theorem~\ref{thm-normalization}]
Let $A\in\Fad$ and put $a_i=q(A_i)$. Let $c=(c_1,\ldots,c_n)$
be the tuple constructed in Section~\ref{sec-normalization}. Since
$q$ is onto, choose $C_i\in\BH$ such that
\begin{align}\label{def-lift-endpoint}
 q(C_i)=c_i,\qquad 1\leq i\leq n.
\end{align}
Put $C=(C_1,\ldots,C_n)$.
The norm-continuous segment $A_i(t)=(1-t)A_i+tC_i$ satisfies
\[
 q(A_i(t))=(1-t)a_i+tc_i=a_i(t).
\]
Proposition~\ref{prop-path} shows that $a(t)$ is doubly commuting and
Taylor invertible. In particular, for $i\ne j$,
\[
 q([A_i(t),A_j(t)])=0,
 \qquad q([A_i(t),A_j(t)^*])=0.
\]
Thus $A(t)\in\Fad$. The uniform estimate follows from
\eqref{eq-uniform-path-bound}, with $\eta=m^2\delta/n$.

Proposition~\ref{prop-endpoint} shows that each $c_i$ is normal and
that $\sum_i|c_i|=1$.  Thus $[C_i,C_i^*]$ is compact for every $i$,
so $C\in\Fen$, and $\sum_i|q(C_i)|=1$.
\end{proof}

To prove Theorem~\ref{thm-components}, we recall the the following theorem of Douglas~\cite{Douglas}, in the
formulation proved in \cite[Theorem 7]{Curto}.

\begin{theorem}\label{thm-douglas-curto}
If $C,D\in\Fen$ and $\ind C=\ind D$, then $C$ and $D$ are joined
by a norm-continuous path in $\Fen$.
\end{theorem}

\begin{proof}[Proof of the equivalence in Theorem~\ref{thm-components}]
The index is constant on norm-continuous paths in $\Fad$, as shown
in Section~\ref{sec-koszul}. Thus condition~(2) implies condition~(1).

Conversely, suppose that $\ind A=\ind B$.
Theorem~\ref{thm-normalization} gives line segments from $A$ to
$C\in\Fen$ and from $B$ to $D\in\Fen$, both contained in $\Fad$.
Index invariance along these segments gives
\[
 \ind C=\ind A=\ind B=\ind D.
\]
By Theorem~\ref{thm-douglas-curto}, $C$ and $D$ can be joined in
$\Fen$.  Concatenating the segment from $A$ to $C$, this path from
$C$ to $D$, and the reversed segment from $D$ to $B$ gives the
required path in $\Fad$.
\end{proof}

Curto \cite[Theorem 3]{Curto} proved surjectivity of the index map
using shifts on the Hardy space of the polydisc. Since these models
are doubly commuting, his examples also give surjectivity on $\Fad$. For completeness,
we include the calculation for the basic shift models and use
finite direct sums to obtain all integer indices.

\begin{proposition}\label{prop-integers}
For every finite $n\geq1$, the index map on $\Fad$ is onto $\Z$.
\end{proposition}

\begin{proof}
Let $\N_0=\{0,1,2,\ldots\}$ and
$H_0=\ell^2(\N_0^n)$, with orthonormal basis
$\{v_{\boldsymbol m}:\boldsymbol m\in\N_0^n\}$.
For the $i$th coordinate vector $\boldsymbol e_i$, define
\[
 S_iv_{\boldsymbol m}=v_{\boldsymbol m+\boldsymbol e_i}.
\]
The operators $S_i$ act in distinct factors of
$\ell^2(\N_0)^{\otimes n}$, so $S=(S_1,\ldots,S_n)$ is doubly
commuting.  Let $P_i$ be the orthogonal projection onto the closed
span of the vectors with $m_i=0$.  Then
\[
 S_i^*S_i=1,\qquad S_iS_i^*=1-P_i.
\]
The projections $P_i$ commute, and their product
$P=P_1\cdots P_n$ is the rank-one projection onto
$\C v_{\boldsymbol0}$. 

Let
$R=\sum_i(1-P_i)$. The vector $v_{\boldsymbol m}$ is an eigenvector
with eigenvalue equal to the number of nonzero coordinates of
$\boldsymbol m$.  It follows that
\begin{align}\label{shift-gap}
 R\geq1-P,\qquad \Ker R=\C v_{\boldsymbol0}.
\end{align}
For $I\subseteq\{1,\ldots,n\}$, the block of $D_S^2$ on the
exterior component indexed by $e_I$ is
\begin{align}\label{shift-block}
 L_I(S)
 &=\sum_{i\in I}S_i^*S_i+\sum_{i\notin I}S_iS_i^*\notag\\
 &=|I|1+\sum_{i\notin I}(1-P_i).
\end{align}
Since $L_I(S)-R=\sum_{i\in I}P_i\geq0$, every block is at least
$R\geq1-P$. Since $P$ has rank one,
the image of every block in the Calkin algebra is at least $1$.
Thus $S$ is Fredholm by Proposition~\ref{prop-laplacian}.
When $I\ne\varnothing$, this block is at least $1$ and has zero
kernel.  When $I=\varnothing$, it equals $R$ and has kernel
$\C v_{\boldsymbol0}$ by \eqref{shift-gap}.  Self-adjointness gives
$\Ker D_S=\Ker D_S^2$, and hence
\[
 \Ker D_S=\C(v_{\boldsymbol0}\otimes1).
\]
Thus, this kernel is entirely in even degree. Define
$B_S=D_S^{\mathrm{odd}}$, we have
\[
 \dim\Ker B_S=0,\qquad \dim\Ker B_S^*=1,
 \qquad \ind S=-1.
\]

Now let $T=(S_1^*,S_2,\ldots,S_n)$, which is again doubly commuting.
Replacing $S_1$ by $S_1^*$ interchanges the two first-coordinate
terms in~\eqref{shift-block}. Thus
\[
 L_I(T)=L_{I\mathbin{\triangle}\{1\}}(S),
\]
where \begin{align}
I\mathbin{\triangle}\{1\}
&=
\begin{cases}
I\cup\{1\},      & \text{if } 1\notin I,\\
I\setminus\{1\}, & \text{if } 1\in I.
\end{cases}
\end{align} The blocks for $T$ are
therefore those for $S$ with their labels changed. In particular,
$T$ is Fredholm, and its only block with nonzero kernel is indexed
by $I=\{1\}$. Hence
\[
 \Ker D_T=\C(v_{\boldsymbol0}\otimes e_1).
\]
This kernel lies in odd degree, so $\ind T=1$.

For a finite direct sum of tuples $A^{(1)},\ldots,A^{(r)}$, the
natural unitary between the exterior tensor spaces identifies
\[
 B_{A^{(1)}\oplus\cdots\oplus A^{(r)}}
 \quad\hbox{with}\quad B_{A^{(1)}}\oplus\cdots\oplus B_{A^{(r)}}.
\]
Hence this operator is Fredholm when each summand is Fredholm, and
its index is the sum of their indices. For $k>0$, take $k$ copies
of $T$; for $k<0$, take $|k|$ copies of $S$. The resulting tuples
are doubly commuting and have index $k$. Their Hilbert spaces are
separable and infinite dimensional, so a unitary identification with
$H$ gives the required tuples on $H$. For $k=0$, use
$(1,0,\ldots,0)$, whose Koszul operator satisfies $D^2=1$.
\end{proof}

This completes the proof of Theorem~\ref{thm-components}.

\begin{corollary}\label{cor-zero}
A tuple $A\in\Fad$ has index zero if and only if it can be joined to
$(1,0,\ldots,0)$ by a norm-continuous path in $\Fad$.
\end{corollary}


\end{document}